\documentclass[12pt]{article}
\usepackage[centertags]{amsmath}     
\usepackage{amssymb}
\usepackage{amsthm}
\usepackage{verbatim}                
\usepackage[dvips]{graphics}
\usepackage{enumerate}

\newtheorem{proposition}{Proposition}
\newtheorem{lemma}[proposition]{Lemma}
\newtheorem{definition}{Definition}

\newtheorem{corollary}[proposition]{Corollary}
\newtheorem{theorem}[proposition]{Theorem}
\newtheorem*{remark}{Remark}

\newcommand{\me}{\mathrm{e}}
\newcommand{\mi}{\mathrm{i}}
\newcommand{\dif}{\mathrm{d}}
\newcommand{\eps}{\varepsilon}

\newcommand{\ud}{u_\mathbb{D}}
\newcommand{\uds}{u_{\mathbb{D}^*}}

\newcommand{\CC}{\mathbb{C}}
\newcommand{\DD}{\mathbb{D}}
\newcommand{\HH}{\mathbb{H}}

\newcommand{\RR}{\mathbb{R}}

\newcommand{\length}{\mathrm{length}}
\newcommand{\dist}{\mathrm{dist}}
\newcommand{\area}{\mathrm{area}}

\begin{document}
\title{Conformal Extension of Metrics of Negative Curvature}
\author{Dan Mangoubi}
\date{\textit {In memory of Robert Brooks}}
\maketitle
\begin{abstract}
We consider the problem of extending a conformal metric
of negative curvature, given outside of a neighbourhood of $0$ in 
the unit disk $\DD$,
to a conformal metric of negative curvature in $\DD$.
We give conditions under which such an extension
is possible, and also give obstructions to such an extension.
The methods we use are based on a maximum principle
and the Ahlfors--Schwarz Lemma.
We also give an example in which no extension is possible,
even when the conformality condition is dropped.
We apply these considerations to compactification of Riemann surfaces.
\end{abstract}

\section{Introduction}
In~\cite{ps}, Brooks considers
the following problem: Let $S^O$ denote 
a complete hyperbolic Riemann surface of finite area,
and let $S^C$ be its conformal compactification. In general,
$S^C$ need not carry a hyperbolic metric, but it is
shown in~\cite{ps} that if the cusps are large,
in a sense to be defined below, 
then the hyperbolic metric on
$S^C$ is close to the 
hyperbolic metric on $S^O$ outside cusp neighbourhoods.
The result in~\cite{ps} is qualitative, and 
in this paper we quantify this result. 
%

First, we would like to know how large a cusp should be
in order to adjust the metric in its neighbourhood
to a non-cusped metric with negative curvature.
For example, if $S^O$ is the complete hyperbolic 
3-punctured sphere, then $S^C$ has no hyperbolic structure
at all. By contrast, we will show in Theorem~\ref{horo-2pi}
below that if a cusp has length $>2\pi$, then we may
modify the metric conformally inside the cusp neighbourhood
to obtain a metric of negative curvature. We will show that the
constant $2\pi$ is sharp.

The above question leads us to the interesting problem 
of extension of metrics, with negative curvature
and with a fixed conformal type, across a given boundary curve.
In this paper, we consider this problem only for
the complete hyperbolic metric on the punctured unit
disk, since our motivation comes from hyperbolic Riemann
surfaces of finite area. There is an obvious necessary
condition for the extension which comes from the Gauss--Bonnet
Theorem (cf.\ \S\ref{sec:question}).
We give several examples which show that this
is not a sufficient condition for 
extension across general curves. Even under convexity
restrictions, which seem natural, the Gauss--Bonnet condition
is not enough.

  A very similar problem was 
already addressed by M.~Gromov
in~\cite{gromov-pdr}, pp.~109--110. Gromov gives an obstruction
for extension of metrics, not necessarily conformal, with curvature
of a fixed sign, which goes beyond the Gauss--Bonnet obstruction.
We present Gromov's obstruction in section~\ref{sec:gromov}.
It seems that in general, it is not easy to decide
whether this obstruction is met.
  We give a class of curves which do satisfy the Gauss--Bonnet necessary 
condition for extension, but meet Gromov's obstruction.

  Next, we estimate how large a cusp should be in order
to adjust the metric in its neighbourhood
to a non-cusped metric with curvature close to $-1$.
This gives a quantitative version of the theorem
from~\cite{ps} which compares between the metrics on $S^O$ and $S^C$.
We will show in Theorem~\ref{crv_estm} that if a cusp has
length of the order $1/\sqrt{\eps}$, then it is possible
to modify the metric in a neighbourhood of the cusp to obtain
a non-cusped metric with Gaussian curvature $\kappa$,  
$$-(1+\eps)<\kappa<-\frac{1}{1+\eps}.$$ We will show that
this estimate is sharp.

\paragraph{Acknowledgements.}
We would like to thank Mikhail Katz for his fruitful ideas.
The main results of this paper were contained in
the author's M.~Sc.\ thesis~\cite{thesis}, written under
the direction of Robert Brooks, whom we owe special thanks
for making several helpful suggestions and clarifications.

I dedicate this paper to the memory of Robert Brooks who 
inspired me greatly, and provided me with enormous 
support and encouragement. I will miss him dearly.

\section{General Settings and the Main Theorems}
Let $\DD^*$ be the complete hyperbolic punctured unit disk
$\{z\in \CC:\, 0<|z|<1\}$ with the unique complete
hyperbolic metric $\dif s^*$, conformally equivalent
to the Euclidean metric. We will give this metric explicitly
later. Let $C_r$ be the Euclidean circle with center $0$ and 
radius $0<r<1$. $C_r$ is a horocycle in $\DD^*$.
Let $B_r=\{z:\, 0<|z|<r\}$ be the interior of $C_r$.

We are interested in adjusting the metric conformally in $B_r$
in such a way that we get rid of the cusp and retain
negative curvature. Thus, we will have a new metric
of negative curvature on the unit disk $\DD$, which coincides with
the metric on $\DD^*$ outside $B_r$.

\begin{definition}
  Let $R$ be a region in\/ $\DD$ which contains $0$.
  We say that a metric $\dif s$ on\/ $\DD$ is a 
  \emph{negatively curved completion} of\/ $\dif s^*$ inside $R$
  if 
  \begin{enumerate}[(i)]
    \item \label{con-class} $\dif s$ is conformally equivalent to\/ $\dif s^*$.
    \item $\dif s = \dif s^*$ outside $R$.
    \item The curvature $\kappa(\dif s)$ of the metric $\dif s$
          satisfies
         $$\kappa(\dif s)< 0.$$
  \end{enumerate}
We may also say that $\dif s$ is \emph{a completion of
$\dif s^*$ across $\partial R$}. 
\end{definition}

The main theorems we prove are
\begin{theorem}
\label{horo-2pi}
$\dif s^*$ has a negatively curved completion inside $B_r$
if and only if $\length(C_r) > 2\pi$.
\end{theorem}

In~\S\ref{sec:question} we give several examples
which show that the horocycle condition in the
above theorem is sharp.
%

The next theorem is in the same spirit, but it gives
a more accurate estimation of the curvature of the adjusted metric.
This theorem leads us in~\S\ref{sec:compact} to a quantitative 
version of the comparison theorem from~\cite{ps}.
\begin{theorem}
\label{crv_estm} There exist\/ $C_1>0$, $C_2>0$ such that,
  \begin{enumerate}[(i)]
    \item  For all\/ $0<\eps<1$, if $r>1-C_1\sqrt{\eps}$, 
           then\/ $\dif s^*$ has a conformal completion\/ 
           $\dif s_\eps$ in $B_r$ with curvature\/
           $$-(1+\eps)<\kappa(\dif s_\eps)<-\frac{1}{1+\eps}.$$
    \item  For all small\/ $\eps>0$, if\ $\dif s^*$ has a conformal
           completion\/ $\dif s_\eps$ in\/ $B_r$ with curvature\/
           $$\kappa(\dif s_\eps)<-\frac{1}{1+\eps},$$
           then\/ $r>1-C_2\sqrt{\eps}$.
  \end{enumerate}
\end{theorem}

\section{Proof of Necessity in Theorem~\ref{horo-2pi}}
In this section, we prove the necessity part of
Theorem~\ref{horo-2pi}. The necessity of the
horocycle condition is a direct consequence of 
the Gauss--Bonnet Theorem.
We denote by $\kappa$ the Gaussian curvature,
and by $\kappa_g$ the geodesic curvature of a curve.
\begin{lemma}
\label{GB-bound}
Let $R$ be a region in\/ $\DD$ which contains\/ $0$.
If\/ $\dif s$ is a negatively curved completion of\/ $\dif s^*$ in $R$,
then 
  $$ \oint_{\partial R} \kappa_g\,\dif(\length) > 2\pi.$$
\end{lemma}
\begin{proof}
  By the Gauss--Bonnet Theorem in the region $R$
equipped with the metric $\dif s$ we have
  $$\iint_{R} \kappa(\dif s) \,\dif(\area) +
    \oint_{\partial R} \kappa_g\, \dif(\length)= 2\pi,$$
and the first term is negative.
\end{proof}
\begin{proof}[Proof of necessity in Theorem~\ref{horo-2pi}]
Let $\dif s$ be a negatively curved completion of $\dif s^*$
in $B_r$.
By the previous lemma 
  \begin{equation}
    \oint_{C_r} \kappa_g\,\dif(\length) > 2\pi.
  \end{equation}
Since the geodesic curvature of a horocycle is $1$,
we obtain $\length(C_r)>2\pi$.
\end{proof}

\section{The Complete Hyperbolic Punctured Disk}
Before proving the sufficiency part in Theorem~\ref{horo-2pi},
we devote the next two sections to present some elementary formulas
in the complete hyperbolic punctured unit disk.

Let $\DD^*$ denote the complete hyperbolic punctured unit disk,
and let $\HH^2$ denote the upper-half-plane $\{z:\,\Im(z)>0\}$
with the hyperbolic metric 
$$\dif s^2=\frac{\dif x^2+\dif y^2}{y^2}.$$
$\DD^*$ is isometric to $\HH^2/\{z\mapsto z+1\}$
via the map 
\begin{equation}
\label{h2-pd}
  z\in\HH^2\mapsto \me^{2\pi\mi z}\in\DD^*.
\end{equation}
This isometry lets us compute easily the metric on $\DD^*$:
\begin{equation}
\label{dstar_metric}
\dif s^*= -\frac{1}{r\log r}|\dif z|,
\end{equation}
where $r=|z|$.

\begin{lemma}
\label{lem:geodesic_curvature}
Let $\gamma$ be a simple closed $C^1$-curve around $0$ in $\DD^*$.
Denote by $\mathrm{int}(\gamma)$ the finite-area component of\/
$\DD^*\setminus\gamma$. The total geodesic curvature
of $\gamma$ is given by
  $$\oint_\gamma \kappa_g\,\dif(\length) = \area(\mathrm{int}(\gamma)).$$
\end{lemma}
\begin{proof}
  The Euler characteristic of $\mathrm{int}(\gamma)$,
   $\chi(\mathrm{int}(\gamma))$, is $0$.
  So, by the Gauss--Bonnet Theorem we have:
  $$\oint_\gamma \kappa_g\,\dif(\length) +
    \iint_{\mathrm{int}(\gamma)} 
    \kappa\,\dif(\area) = 2\pi\chi(\mathrm{int}(\gamma))=0.$$
  The lemma follows at once since $\kappa=-1$.
\end{proof}

\begin{lemma}
  \label{lem:area-per}
  $\area(B_r)= \length(C_r) = -2\pi/\log r$.
\end{lemma}
\begin{proof}
By~(\ref{dstar_metric}),
$$
  \begin{array}{lclclcl}
    \area(B_r) &=&
      \displaystyle\int_0^{2\pi}\!\!\!\int_0^r \frac{t}{(t\log t)^2}
                    \:\dif t\:\dif\theta
        &=& -2\pi\left.\displaystyle\frac{1}{\log t}\right|_0^r
        &=& \displaystyle{-\frac{2\pi}{\log r}.}
  \end{array}
$$
Next, we can compute the length of $C_r$ directly, or we can apply
the Gauss--Bonnet Theorem to show that it is equal to the area
(recall that $\kappa_g=1$ for $C_r$).
\end{proof}

\begin{figure}
  \begin{center}
  \scalebox{.6}{\includegraphics{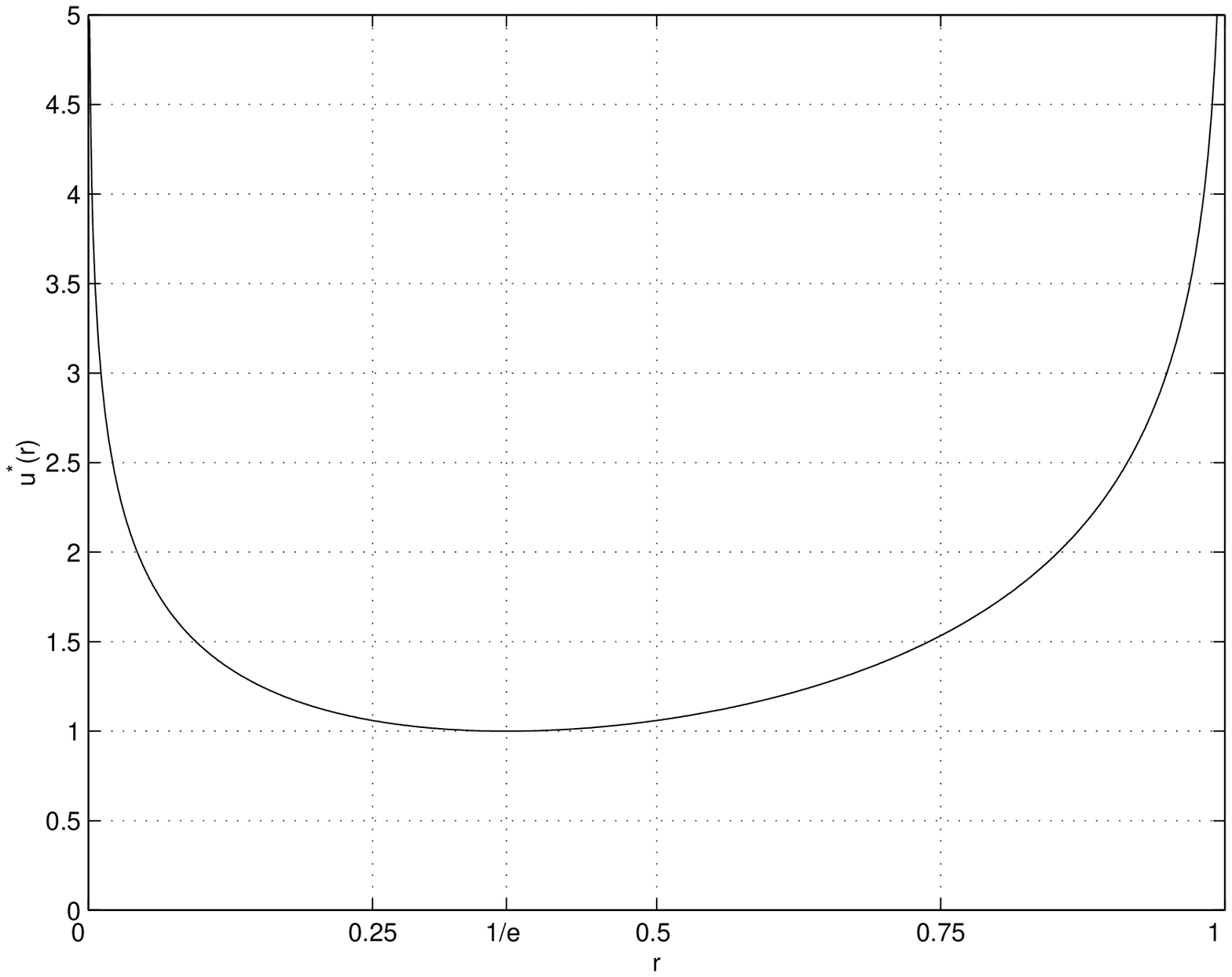}}
  \end{center}
  \caption{The function $u^*$, where $\dif s^*=\me^{u^*}|\dif z|$}
  \label{fig:udstar}
\end{figure}

\begin{figure}
  \begin{center}
  \scalebox{.6}{\includegraphics{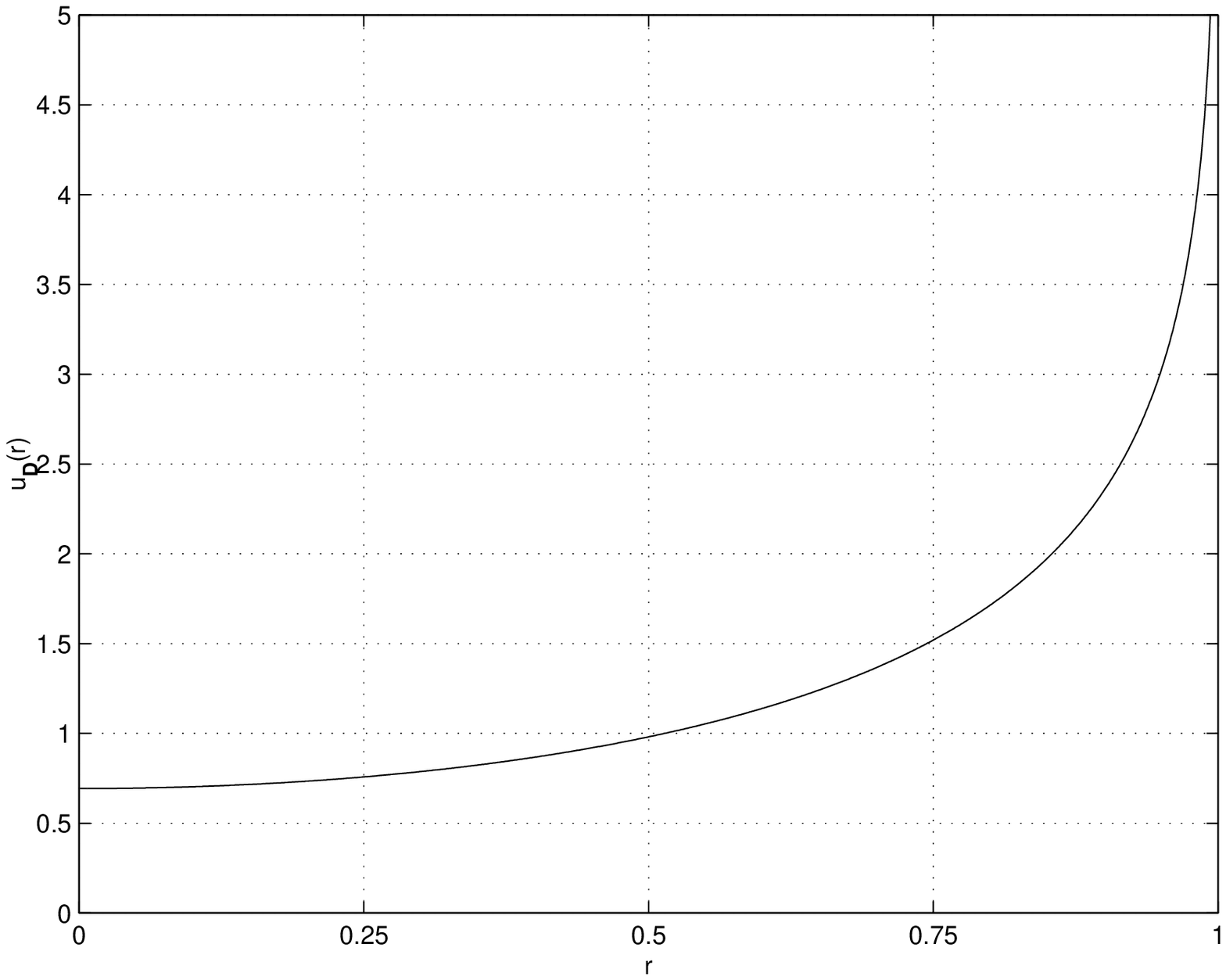}}
  \end{center}
  \caption{The function $u_\DD$, where $\dif s_\DD=\me^{u_\DD}|\dif z|$}
  \label{fig:ud}
\end{figure}

We would like now to see some simple properties
of the metric on $\DD^{*}$ given in~(\ref{dstar_metric}).
Set
\begin{eqnarray}
  \lambda^*(r) &=& -\frac{1}{r\log r},\\
  u^*(r) &=& \log(\lambda^*(r))=\log\left(\frac{1}{r\log(\frac{1}{r})}\right).
  \label{eqn:udstar}
\end{eqnarray}
The function $u^*$ is shown in Figure~\ref{fig:udstar}.
For comparison, we show in Figure~\ref{fig:ud}
the corresponding function $u_\DD$
for the complete hyperbolic metric $\dif s_{\DD}$ on the unit disk:
  $$u_\DD = \log\left(\frac{2}{1-r^2}\right).$$
\begin{lemma}
  \label{lem:udstar}
  $u^*$ satisfies 
  \begin{enumerate}[(i)]
    \item $\lim_{r\to 0} u^*(r) =\infty.$
    \item $\lim_{r\to 1} u^*(r) =\infty.$
    \item $(u^*)'(r)<0$ for $r<\frac{1}{\me}$.
    \item $(u^*)'(r)>0$ for $r>\frac{1}{\me}$.
    \item For all $0<r<1$, $(u^*)''(r)>0$.
  \end{enumerate}
\end{lemma}
\begin{proof}
The proof is a simple computation.
\end{proof}

\section{Curvature for Conformally Euclidean Metrics}
We recall that if 
  $$\dif s^2 = \lambda(x,y)^2(\dif x^2 +\dif y^2) $$
is a conformally Euclidean Riemannian metric in a 
domain $U\subseteq \RR^2$,
then its curvature is given by the well known formula
\begin{equation}
  \label{curvature}
  \kappa(\dif s)=-\frac{\Delta\log\lambda}{\lambda^2},
\end{equation}
where $\Delta=\partial^2/\partial x^2+\partial^2/\partial y^2$.

In particular, for a radially symmetric Riemannian metric 
  $$\dif s=\me^{u(r)}|\dif z|$$
on the unit disk $\DD$,~(\ref{curvature})~reduces to 
\begin{equation}
  \label{curv}
  \kappa(\dif s)=-\frac{u''+\frac{1}{r}u'}{\me^{2u}}.
\end{equation}

\section{Proof of Sufficiency in Theorem~\ref{horo-2pi}}
\label{sec:suff}
We are ready to finish the proof of Theorem~\ref{horo-2pi}.
\begin{proof}[Proof of sufficiency in Theorem~\ref{horo-2pi}]
Suppose $\length(C_{r_0})>2\pi$. 
By Lemma~\ref{lem:area-per}
$r_0>\frac{1}{\me}$. Then, Lemma~\ref{lem:udstar} tells us that
$(u^*)'(r_0)>0$ and $(u^*)''(r_0)>0$.  
Therefore, there exists a convex monotonically increasing
$C^2$-function 
$u(r)$ on $[0,1)$ such that
$$ 
 \begin{array}{l}
  u(r)=u^*(r) \quad\mbox{for}\ r>r_0,\\
  u'(0)=0.
\end{array}
$$
Indeed, we can construct $u(r)$ as follows:
First, define a positive continuous function $w(r)$ such that,
\begin{eqnarray*}
  \int_0^{r_0} w(r)\,\dif r & = & (u^*)'(r_0), \\ 
  w(r)  & = & (u^*)''(r)\quad\mbox{for}\ r\geq r_0.
\end{eqnarray*}
Next, define the function $v(r)$ by
\begin{equation*}
  v(r)=\int_0^r w(s)\,\dif s. 
\end{equation*}
Observe that $v(r)$ is a $C^1$-function which coincides with $(u^*)'(r)$
for $r\geq r_0$.
Finally, define $u(r)$ by
\begin{equation*}
  u(r)={u^*}(r_0)+\int_{r_0}^r v(s)\,\dif s.
\end{equation*}
$u$ is a $C^2$-function which satisfies:
\begin{eqnarray*}
  \forall r\geq r_0,\ u(r) &=& {u^*}(r),\\
  u'(0)&=&v(0)=0,\\
  \forall r>0,\ u'(r)>0 &\Rightarrow&\mbox{$u$ is monotonically increasing},\\
  u''=w>0 &\Rightarrow& \mbox{$u$ is convex}.
\end{eqnarray*}

To conclude, set $\lambda(r)=\me^{u(r)}$, and
$\dif s=\lambda(r)|\dif z|$.
Then, we have:
\begin{enumerate}[a)]
  \item $\lambda'(0)=0 \Rightarrow$ $\dif s$ is smooth at $0$.
  \item $\dif s$ coincides with $\dif s^*$ for $r>r_0$.
  \item $\kappa(\dif s)<0$ since $\log\lambda = u$ is convex,
        monotonically increasing, and since the curvature 
        is given by~(\ref{curv}).
\end{enumerate}
\end{proof}
\section{A Natural Question}
\label{sec:question}
Lemma~\ref{GB-bound} raises the possibility that there might be
a negatively curved completion of $\dif s^*$
in any region $R$\/ in $\DD$ which contains $0$ and satisfies
the Gauss--Bonnet restriction
\begin{equation}
  \label{GB-res}
  \oint_{\partial R} \kappa_g\,\dif(\length) > 2\pi,
\end{equation}
which would generalize Theorem~\ref{horo-2pi}.
In this section we show, that the supposedly more
general theorem is false. We build counterexamples
of several different types.
\subsection{A First Example}
We give an example which shows that 
the Gauss--Bonnet restriction~(\ref{GB-res}) is not enough 
for having a negatively curved completion of $\dif s^*$ in $R$. 
We use the following general criterion for not having a completion:
\begin{theorem}
\label{criterion}
Let $\gamma$ be a curve around\/ $0$ in\/ $\DD^*$ such that
the function\/ $\lambda^*(z)$ attains its maximum on $\gamma$
at a point $z_0$ for which $|z_0|<1/\me$.
Then, $\dif s^*$ has no negatively curved completion
across $\gamma$.
\end{theorem}
\begin{proof}
  Suppose, on the contrary, that
  $\dif s = \lambda(z)|\dif z|$
  is a negatively curved completion of $\dif s^*$
  across $\gamma$.
    On the one hand, since $\kappa(\dif s)<0$, 
  formula~(\ref{curvature}) tells us
  that $u=\log\lambda$ is a subharmonic function.
  Hence, by the maximum principle, $u(z_0)$
  is the maximal value of $u$
  in $\overline{\mathrm{int}(\gamma)}$.
  So, if $\hat{n}$ is the outer unit normal to $\gamma$, then
  \begin{equation}
   \label{nor_der}
   \frac{\partial u}{\partial \hat{n}}(z_0)\geq 0.
  \end{equation}
  On the other hand by Lemma~\ref{lem:udstar}, 
   \begin{equation}
     \label{r_der}
     \frac{\partial u}{\partial \hat{r}}(z_0) < 0,
   \end{equation}
   where $\hat{r}$ is a unit vector from $0$ to $z_0$.
   Since $z_0$ is a maximal point, we have
    \begin{equation}
      \label{r-n}
      \frac{\partial u}{\partial \hat{r}}(z_0) = 
      \frac{\partial u}{\partial \hat{n}}(z_0)\langle\hat{r}, \hat{n}\rangle,
    \end{equation}
  Now, notice that by Lemma~\ref{lem:udstar}, the ray from $0$ to $z_0$ 
  is contained in $\mathrm{int}(\gamma)$. Therefore, 
  $\langle\hat{r}, \hat{n}\rangle>0$. 
  From here we get that~(\ref{nor_der}) and~(\ref{r-n})
  are in contradiction with~(\ref{r_der}).
\end{proof}

\begin{theorem}
\label{simple_ex}
There exists a simple closed smooth curve $\gamma$, homologous
to $C_r$\/ in $\DD^*$, such that
\begin{enumerate}[(i)]
  \item \label{nec} $\oint_{\gamma} \kappa_g\,\dif(\length) > 2\pi.$
  \item $\dif s^*$ has no negatively curved completion across $\gamma$.
\end{enumerate}
\end{theorem}
\begin{proof}
  \begin{figure}
    $$\includegraphics{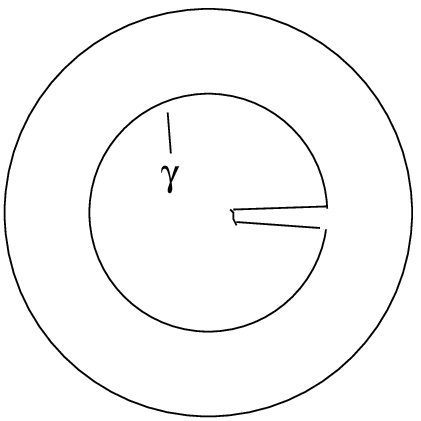}$$
    \caption{horocycle with a slit}
    \label{fig:horocycle_slit}
  \end{figure}
Let $\gamma$ be a horocycle with a slit
(see Figure~\ref{fig:horocycle_slit}), i.e.
smooth the curve which is the boundary of the region
  $$ V=\left\{|z|<R_2\right\}\setminus
   \left\{|z|>R_1 \mbox{ and } |\arg z|<\theta\right\},$$
  with $R_1,$ $R_2$, $\theta$ chosen as follows:
  Let $R_2$ satisfy $\frac{1}{\me}<R_2<1$. 
  By Lemma~\ref{lem:udstar} we can find $0<R_1<\frac{1}{\me}$ such that 
  \begin{equation}
    u^*(R_1)>u^*(R_2).
  \end{equation}
  By Lemma~\ref{lem:area-per}
  $\area(B_{R_2}) > 2\pi$. Therefore,
  we can find $\theta$ small enough, such that the area of $V$
  $>2\pi$.
  
  Lemma~\ref{lem:geodesic_curvature} assures us that the Gauss--Bonnet 
  restriction~(\ref{nec}) is satisfied, 
  and from Theorem~\ref{criterion}, we immediately get that
  $\dif s^*$ has no negatively curved completion across $\gamma$.
\end{proof}

\subsection{A Convex Example}
\label{sec:convex}
Now, we impose convexity on our curve.
A convex curve is a curve whose geodesic curvature is
non-negative.

\begin{theorem}
There exists a convex simple closed smooth curve $\gamma$, homologous
to $C_r$ in\/ $\DD^*$, such that
\begin{enumerate}[(i)]
  \item $\oint_{\gamma} \kappa_g\,\dif(\length) > 2\pi.$
  \item $\dif s^*$ has no negatively curved completion across $\gamma$.
\end{enumerate}
Moreover, such a curve can be taken arbitrarily close to the horocycle
of length $2\pi$.
\end{theorem}
\begin{proof}
  Take $\gamma$ to be composed of a horocyclic segment and 
  a geodesic segment
  as in Figure~\ref{fig:curve}:
  \begin{figure}
    $$\scalebox{.75}[.75]{\includegraphics{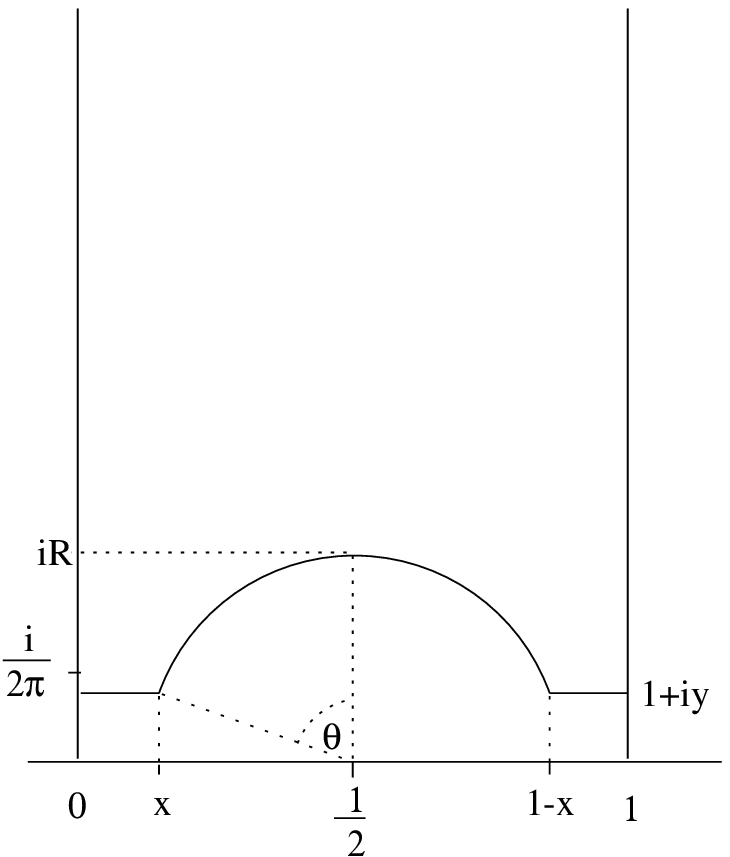}}$$
    \caption{A convex curve}
    \label{fig:curve}
  \end{figure}
  $\gamma$ is the boundary of the domain
  $$ V = \left\{ z\in \mathbb{H}^2 :\, |z -\frac{1}{2}|^2 > R^2,
          \ 0\leq\Re z < 1,\ \Im z > R\cos\theta \right\} $$
  The total geodesic curvature of $\gamma$ is given by the area of $V$ 
  (see Lemma~\ref{lem:geodesic_curvature}):
  $$\oint_\gamma\kappa_g\,\dif\mbox{(length)} = 
    \frac{1-2R\sin\theta}{R\cos\theta}+2\theta,$$
  which we would like to be $>2\pi$ (for convenience, we write
  this in terms of $y=R\cos\theta$):
  \begin{equation}
    \label{ineq:k_g>2pi}
    \frac{1}{y}-2\tan\theta+2\theta>2\pi.
  \end{equation}
  In order to apply the non-extendibility criterion of
  Theorem~\ref{criterion},
  we should also like to have (recall the map in Equation~(\ref{h2-pd}))
  \begin{eqnarray}
     u^*(\me^{-2\pi R}) &>& u^*(\me^{-2\pi y})\nonumber\\
     &\Leftrightarrow& \nonumber\\
     \label{ineq:max}
     \frac{\me^{2\pi R}}{2\pi R} &>&
     \frac{\me^{2\pi y}}{2\pi y}.
  \end{eqnarray}
  We may first find a pair $y_0$, $\theta_0>0$ which satisfies
  equality in~(\ref{ineq:k_g>2pi}) and satisfies inequality~(\ref{ineq:max}),
  and then decrease $\theta_0$ a little, keeping $y_0$ fixed.
  So, extracting $y$ from~(\ref{ineq:k_g>2pi}), we obtain that we
  should find $0<\theta<\pi/2$ such that
  \begin{equation}
    \frac{\pi(1-\cos\theta)}{\pi+\tan\theta-\theta}>
    -\cos\theta\,\log\cos\theta.
  \end{equation}
  The existence of such a $\theta$ can be shown by elementary calculus.
  In fact, the first three derivatives of the difference are $0$ 
  at $\theta=0$, and the fourth derivative of 
  the difference at $\theta=0$ is $>0$.
  So, we can find $\theta$ arbitrarily close to $0$, which satisfies
  the last inequality. This shows that we can take our
  curve arbitrarily close to the horocycle of length $2\pi$, as stated.
\end{proof}

\subsection{The Condition of Gromov}
\label{sec:gromov}
In \cite{gromov-pdr}, pp.\ 109--110, Gromov gives a
necessary condition for extension of metrics, which is 
based on the Gauss--Bonnet Theorem. The condition as stated
in Gromov's book is a little obscure, and we rephrase it
here for the sake of clarity:
\begin{theorem}[\cite{gromov-pdr}, pp.\ 109--110]
Let $\gamma$ be a simple closed smooth curve bounding a disk $V$.
Suppose that outside of $V$, we are given a Riemannian
metric $\dif s$ of negative curvature.
Let $p, q\in \gamma$. $p, q$ dissect $\gamma$ into two segments
$\alpha$ and $\beta$. Let $\length(\beta)\le\length(\alpha)$.
 Suppose that every $\gamma^+\subseteq\alpha$, possibly
disconnected, with $\length(\gamma^+)\le\length(\beta)$,
satisfies
  $$\int_{\alpha\setminus\gamma^+} \kappa_g\,\dif(\length)\le 0.$$
Then, $\dif s$ cannot be extended to a 
negatively curved metric in $V$.
\end{theorem}
\begin{proof}
Let $\dif s$ be a negatively-curved smooth Riemannian metric in $V$, which
extends the given metric outside of $V$.
Let $\delta$ be a path of minimal length 
in $\overline{V}$ from $p$ to $q$ which is different from $\alpha$.
Then,
\begin{enumerate}
\item\label{smooth}
 $\delta$ is a $C^1$ curve,
\item\label{kg<0}
 $\delta\cap\partial V$ has non-positive geodesic curvature.
\end{enumerate}
Define $\gamma^+=\delta\cap\alpha$.
Apply the Gauss--Bonnet Theorem to each of the complementary
regions of $V\setminus\delta$ which touch $\alpha\setminus\gamma^+$.
Taking negative curvature into consideration,
together with~\ref{smooth}.\ and~\ref{kg<0}., we get
  $$\int_{\alpha\setminus\gamma^+}\kappa_g > 0,$$
contradicting our assumptions.
\end{proof}

Gromov also gives obstructions for
extension of metrics of positive curvature,
and for arbitrary Euler characteristic of $V$.
These cases are treated similarly, and we will omit them here.

We give a class of examples which 
satisfy the Gauss--Bonnet necessary condition
for extension, but meet Gromov's obstruction.

The idea, shown to us by Mikhail Katz, is to find
a situation in which two points are connected by
two geodesics. This is in contradiction with
the fact that in simply connected non-positively 
curved spaces there is at most one geodesic 
connecting two points. Alternatively, one can see the following
as a special case of Gromov's obstruction, taking $\alpha$
to be the geodesic.
\begin{theorem}
\label{thm:no-extend_all}
There exists a convex simple closed smooth curve $\gamma$, homologous
to $C_r$ in $\DD^*$, such that
\begin{enumerate}[(i)]
  \item $\oint_{\gamma} \kappa_g\,\dif(\length) > 2\pi.$
  \item $\dif s^*$ has no negatively curved completion across $\gamma$,
        even if we allow completions
        with metrics not conformal to the Euclidean metric.
\end{enumerate}
\end{theorem}
\begin{proof}
We construct $\gamma$ as follows (see Figure~\ref{fig:no_extend_all}):
\begin{figure}
  $$\includegraphics{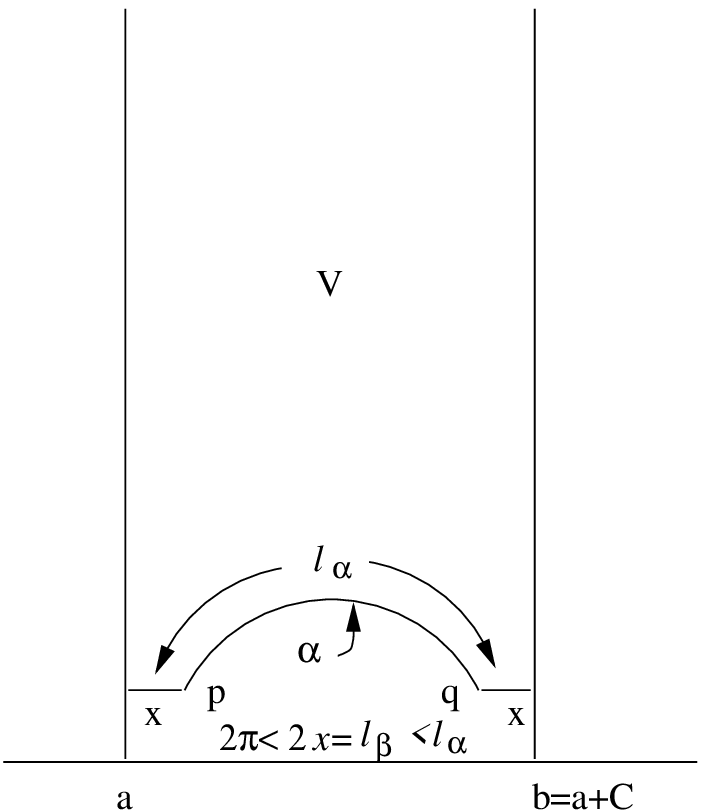}$$
  \caption{No extension}
  \label{fig:no_extend_all}
\end{figure}
Take any geodesic $\tau$ in $\HH^2$, 
whose end points are in $\RR$.
Take a symmetric geodesic segment $\alpha$ about 
the ``center'' of $\tau$ of hyperbolic length
$l_{\alpha}>2\pi$. Now, extend the end points $p, q$ 
of $\alpha$ by horocyclic segments
of total hyperbolic length $2\pi<l_{\beta}<l_{\alpha}$. 
Then, identify $\DD^*$ with 
$\mathbb{H}^2/\{z\mapsto z+C\}$, where $C$ is the Euclidean distance
between the end points of the horizontal segments.
We get a convex curve $\gamma$ in $\DD^*$ 
whose total geodesic curvature $>2\pi$ (Lemma~\ref{lem:geodesic_curvature}).

Suppose $\dif s$ is a negatively curved completion of $\dif s^*$ 
across $\gamma$ (it might not be conformally Euclidean). 
Since the distance between $p$ and $q$
\begin{equation}
\label{two_geo}
  \dist(p, q)\le l_{\beta}<l_{\alpha},
\end{equation}
there exists a geodesic $\delta\neq\alpha$ which connects these points.
But this is in contradiction with the fact that in a Riemannian surface of
negative curvature there is at most one geodesic in each homotopy class.
\end{proof}

\begin{remark}\upshape
The above construction shows the existence of such
curves. We could also work as in the previous example
in the fixed model $\HH^2/\{z\mapsto z+1\}$. In that case, the existence
is less obvious: Fixing $\theta$ (see Figure~\ref{fig:curve})
determines the length of the geodesic segment $\alpha$. Then, 
we can take concentric segments until we get the desired 
inequality~(\ref{two_geo}).
In terms of the last example, we should find $R, \theta$ which
satisfy the inequalities~(\ref{ineq:k_g>2pi}) and
\begin{equation}
  \label{ineq:geo>horo}
  \frac{1}{\cos^2(\theta)}>
  \cosh\left(\frac{1}{R\cos\theta}-2\tan\theta\right).
\end{equation}
Here, the left-hand side is the hyperbolic cosine of the length of the
geodesic segment $\alpha$, 
while the right-hand side is the hyperbolic cosine of
the length of the horocyclic segment.
\end{remark}
%

\subsection{A Remark on the Last Two Examples}
The examples presented in sections~\S\ref{sec:convex}
and~\S\ref{sec:gromov} seem similar, but in fact
there are curves which satisfy the non-extendibility 
criterion of~\S\ref{sec:convex} but not that of~\S\ref{sec:gromov},
and vice versa. There are also curves which satisfy simultaneously
the non-extendibility criteria of~\S\ref{sec:convex} 
and~\S\ref{sec:gromov}.

To see this, we will give values for $R,\ \theta$ which
either satisfy or do not satisfy the 
inequalities~(\ref{ineq:k_g>2pi}),~(\ref{ineq:max})
and~(\ref{ineq:geo>horo}).
The curve with parameters (see Figure~\ref{fig:curve}) 
$R=0.2,\ \theta=0.8$ satisfies only the non-extendibility criterion
of~\S\ref{sec:convex}. The curve with parameters $R=0.48,\ \theta=1.55$
satisfies only the non-extendibility criterion
of~\S\ref{sec:gromov}. Finally, the curve with parameters 
$R=0.4,\ \theta=1.45$ satisfies both of the non-extendibility criteria.

\section{Curvature Estimations}
\label{sec:curvature}
In this section we prove Theorem~\ref{crv_estm}.
Until now, we were interested in adjusting the metric in
$B_r$ to a metric of negative curvature.
Now, we want to have more control on the curvature of
the new metric. Namely, we want it to be close to~$-1$.
We begin by an application of the Ahlfors--Schwarz Lemma~(\cite{as-lemma}).
\subsection{An Ahlfors--Schwarz Bound}
In this section we find a lower bound for $r$ in
order to have a completion of $\dif s^*$ in $B_r$ 
with curvature as in Theorem~\ref{crv_estm}.
The proof is a direct application of the Ahlfors--Schwarz Lemma.
\begin{proof}[Proof of Theorem~\ref{crv_estm}, part(ii)]
Suppose $\dif s_{\eps}$ is a completion 
of $\dif s^*$ in $B_{r_\eps}$ with curvature
  \begin{equation}
    \kappa(\dif s_{\eps}^2) <-\frac{1}{1+\eps}.
  \end{equation}
The last inequality may also be written as
  \begin{equation}
    \kappa(\dif s_{\eps}^2) <\kappa((1+\eps)\dif s_{\DD}^2)<0.
  \end{equation}
Hence, from the Ahlfors--Schwarz Lemma, we obtain
  \begin{equation}
     \frac{\dif s_\eps^2}{\dif s_\DD^2} <1+\eps.
  \end{equation}
For $r>r_\eps$ it means that
  \begin{equation}
     \frac{(1-r^2)^2}{4r^2(\log r)^2}<1+\eps.
  \end{equation}
Expanding the left-hand side near $r=1$ gives
  \begin{equation}
     1+\frac{1}{3}(1-r)^2 + o((1-r)^2) < 1+\eps.
  \end{equation}
For an inequality of this sort to be valid, we must have
  $$r_\eps>1-C\sqrt{\eps}$$
for all small $\eps>0$, with $C>0$ which does not depend on $\eps$. 
\end{proof}
\subsection{Proof of Theorem~\ref{crv_estm}}
Here we construct a completion of $\dif s^*$ in $B_r$,
and we estimate its curvature. The existence of the constant $C_1$
as stated in Theorem~\ref{crv_estm} will follow.

We consider metrics with radial symmetry. Let
    $$\dif s =\lambda(r) |\dif z|$$
be such a metric on $\DD$.
We know that its curvature is given by
\begin{equation}
\label{curv_form}
  \kappa(\dif s) = -\frac{u''+\frac{1}{r}u'}{e^{2u}},
\end{equation}
where $u=\log\lambda$.

  We denote by $\ud$, $\uds$ the $u$'s which correspond
to the complete hyperbolic metrics on $\DD$ and on $\DD^*$,
respectively.
Explicitly (see Figure~\ref{fig:udstar}),
 \begin{eqnarray}
   \ud &=& \log \left(\frac{2}{1-r^2}\right),\\
   \uds &=& \log \left(\frac{1}{r\log(\frac{1}{r})}\right).
 \end{eqnarray}

The idea is to find
an intermediate function $u_\eps$ which coincides with $\uds$
near $r=1$, finite near $r=0$, and is close
to $\ud$ or $\uds$ for any $r$.

We begin with a technical lemma:
\begin{lemma}
\label{taylor}
$\exists A>0$ such that for $r$ sufficiently close to 
$1$ we have:
\begin{enumerate}
  \item $0<\uds(r)-\ud(r)<A(1-r)^2$,
  \item $1-A(1-r)^2<\frac{\uds'(r)}{\ud'(r)}<1$,
  \item $1<\frac{\uds''(r)}{\ud''(r)}<1+A(1-r)^2$.
\end{enumerate}
\end{lemma}
\begin{proof}
  From the Taylor expansions near $r=1$, we see that any $A>1/3$ will do.
\end{proof}
%
%
Let $r_{\eps} = 1-\sqrt{\eps/A}$.
Define the intermediate function $u_\eps$ in two steps. First,
define
\begin{equation}
  \tilde{u}_\eps(r) = \left\{\begin{array}{lcl}
          \uds(r) &, & r_{\eps}\leq r<1\\
          \ud(r) + \uds(r_\eps)-\ud(r_\eps) &, & 0\leq r<r_{\eps}.
          \end{array}\right.
\end{equation}
Then, smooth $\tilde{u}_\eps$\/ near $r(\eps)$\/ without changing much 
its first and second derivatives.
\begin{lemma}
\label{ueps}
The function $u_\eps$ satisfies
\begin{enumerate}
  \item $0<u_\eps-\ud<\eps$,
  \item $1-\eps<\frac{u_\eps'}{\ud'}<1$,
  \item $1<\frac{u_\eps''}{\ud''}<1+\eps$
\end{enumerate}
for all $\eps$ small enough.
\end{lemma}
\begin{proof}
  This follows immediately from the construction of $u_\eps$
and from Lem\-ma~\ref{taylor}.
\end{proof}

We are ready to prove the theorem.
\begin{corollary}
The curvature $\kappa$ of the metric 
  $$\dif s_\eps = \me^{u_{\eps/4}(r)}|\dif z|$$
satisfies 
  $$ -(1+\eps)<\kappa(\dif s_\eps)< -\frac{1}{1+\eps},$$
for all $\eps$ small.
\end{corollary}
\begin{proof}
  We have by formula~(\ref{curv_form})
  $$\frac{\kappa(\dif s_{\eps})}{\kappa(\dif s_\DD)} =
    \frac{u_{\eps/4}''+ \frac{1}{r}u_{\eps/4}'}
    {\ud''+\frac{1}{r}\ud'}\cdot\frac{1}{\me^{2(u_{\eps/4}-\ud)}}.$$
Then, by Lemma~\ref{ueps} the following inequalities are true for 
small $\eps$:
  \begin{equation*}
    \frac{\kappa(\dif s_{\eps})}{\kappa(\dif s_\DD)} >
    \frac{1-\eps/4}{\me^{\eps/2}}>\frac{1}{1+\eps},
  \end{equation*}
  and
  \begin{equation*}
    \frac{\kappa(\dif s_{\eps})}{\kappa(\dif s_\DD)}
         < 1+\eps/4<1+\eps.
  \end{equation*}
Since $\dif s_\eps$ coincides with $\dif s^*$ for $r>1-C\sqrt{\eps}$,
where $C=1/(2\sqrt{A})$ and $A$ is from Lemma~\ref{taylor},
the last two inequalities finish the proof of Theorem~\ref{crv_estm}, part~(i).
\end{proof}

\section{Compactification of Riemann Surfaces}
\label{sec:compact}
Let $S^O$ be a complete hyperbolic Riemann surface of finite area.
We compactify $S^O$ conformally, 
i.e.\ we simply fill in the cusps, 
and take the unique conformal structure on
the filled in surface
to get a compact Riemann surface $S^C$.
We would like to adjust the metric on $S^O$ in disjoint neighbourhoods 
of the cusps to get a smooth metric across the cusps and 
retain negative curvature.
\begin{definition}
We say that the length of a cusp, $P$, on $S^O$, is
$\geq l$, if there is a neighbourhood of $P$,
which is isometric to a domain $V$ in the punctured unit disk,
which contains a closed horocycle of length $\geq l$.
\end{definition}
\begin{definition}
We say that the lenghts of the cusps on $S^O$ is $\geq l$,
if there exist disjoint closed horocycles around each cusp, 
each of length $\geq l$.
\end{definition}

The following theorem is an immediate consequence
of Theorem~\ref{crv_estm}.
\begin{theorem}
Let $S^O$ be a complete hyperbolic finite-area Riemann surface.
We can adjust the metric in horocyclic neighbourhoods
of the cusps to get a smooth metric of negative curvature on $S^C$
if and only if the cusps of $S^O$ have lengths $\geq 2\pi$.
\end{theorem}

In~\cite{ps}, Brooks compares between the complete 
hyperbolic metrics $\dif s^o$, $\dif s^c$ on $S^O$ and $S^C$,
respectively. 
His result is based on the Ahlfors--Schwarz Lemma. The curvature estimates
in~\S\ref{sec:curvature} give a quantitative version of this result:
\begin{theorem}[compare \cite{ps}]
  For every $\eps>0$, there exists $l(\eps)=O(1/\sqrt{\eps})$
  as $\eps\to 0$, such that 
  for any complete hyperbolic finite-area Riemann surface $S^O$
  with cusps $\geq l(\eps)$, we have 
  $$\frac{1}{1+\eps}(\dif s^o)^2\leq(\dif s^c)^2\leq(1+\eps)(\dif s^o)^2, $$
  outside horocyclic neighbourhoods of length $l(\eps)$.
\end{theorem}
%

\bibliographystyle{amsalpha}
\nocite{as-lemma, thesis}
\bibliography{$BIBTEX/phd}
\end{document}